\documentclass[12pt,twoside]{article}

\usepackage[english]{babel}
\usepackage{amsmath}
\usepackage{amsfonts,dsfont}
\usepackage{amssymb}
\usepackage{enumerate}
\usepackage{mathrsfs}
\usepackage{amssymb}
\usepackage[all]{xy}
\usepackage{graphics}
\usepackage{pstricks,pst-plot,pst-node}

\usepackage{dsfont}
\usepackage{amsthm}

\date{}

\begin{document}

\newtheorem{theorem}{\quad theorem}[section]

\newtheorem{Definition}[theorem]{\quad Definition}

\newtheorem{Corollary}[theorem]{\quad Corollary}

\newtheorem{lemma}[theorem]{\quad lemma}

\centerline{}

\centerline{}

\centerline {\bf On the Diameter of the Commuting Graph of a Full Matrix Ring over a Division Ring.}

\centerline{}

\centerline{\bf {C. Miguel}}

\centerline{}

\centerline{ Instituto de Telecomunuca\c
c\~oes,}

\centerline{P\'olo de Covilh\~a,}

\centerline{celino@ubi.pt}

\begin{abstract}For a division ring $\mathds D$, finite dimensional over its center $\mathds F$, we prove that the commuting graph $\Gamma(M_n(\mathds D))$, where $n\geq 3$, is connected
 if and only if the following condition is satisfied

$(\ast)$For any noncentral  matrix $A\in M_n(\mathds D)$ if the algebra $\langle A\rangle$ is a field extension of $\mathds F$ of degree $nk$, then there exists a proper
intermediate field between $\mathds F$ and $\langle A\rangle$.

 Furthermore, if the graph $\Gamma(M_n(\mathds D))$ is connected we prove that its diameter is between four and six. 

\end{abstract}



\section{Introduction }

In a graph $G$, a  path  is an ordered   sequence  $v_1-v_2-\ldots -v_l$ of distinct vertices of $G$  in which every two consecutive vertices are adjacent. The graph $G$ is called  connected if there is at least one path between every pair of its vertices.
 The distance between two distinct vertices $u$ and $v$, denoted by $d(u , v)$, is the length of the shortest  path connecting them (if
such a path does not exists, then $d(u , v)=\infty$). The  diameter of a graph $G$ is the longest distance between any two distinct vertices of  $G$ and will be denoted by $diam(G)$.

The commuting graph of a noncommutative ring $R$, denoted by $\Gamma(R)$, is a graph whose vertices are all noncentral elements of $R$ and two distinct vertices $a$ and $b$ are adjacent if and only if $ab=ba$. In particular, the set of neighbors of a vertex $a$ is the set of all noncentral elements in the centralizer of $a$ in $R$, that is, of $C_R(a)=\{x\in R : ax=xa\}$. When the ring
$R$ is a matrix ring over a field $\mathds F$ the connectedness and diameter of the commuting graph   have
been  studied extensively, see e.g. \cite{ak2, ak3, ak4, da, dol, ga, ce}. If the field $\mathds F$ is algebraically closed, then the commuting graph $\Gamma(M_n(\mathds F))$ is connected and its diameter is always equal to four, provided $n\geq 3$, see \cite{ak2}. Note that for $n=2$ the commuting graph is always disconnected \cite[remark 8]{ak4}.  If the field $\mathds F$ is not algebraically closed, the commuting graph $\Gamma(M_n(\mathds F))$ may be disconnected for an arbitrarily large integer $n$, see \cite{ak3}. However, in $\cite{ak2}$ it is proved that for any field $\mathds F$ and $n\geq 3$, if $\Gamma(M_n(\mathds F))$ is connected, then the diameter is between four and six.

Extensions of results about commuting graphs of matrix rings to the case where the underlying field is replaced by a division ring $\mathds D$ has been of interest, see e.g.
\cite{ak2, ak4}. In particular, it   was determined the diameters of some induced subgraphs of $\Gamma(M_n(\mathds D))$, such as the induced subgraphs on  the set of all  non-scalar non-invertible, nilpotent, idempotent, and involution matrices in $M_n(\mathds D)$.
In this paper we investigate the  connectedness and diameter of the commuting graph $\Gamma(M_n(\mathds D))$ of a full matrix ring over a division ring $\mathds D$, finite dimensional over its center $\mathds F$. We  will establish a necessary and sufficient condition for the connectedness of the graph $\Gamma(M_n(\mathds D))$ and we will prove that if this graph is connected, then its diameter is between four and six, as in the  field case.

\section{Preliminaries}
In this section we assemble the tools that we require to prove your results. We begin with the rational canonical form for matrices over a division ring.
 Remember that if $A$ is a $r\times r$ matrix and $B$ is a $s\times s$ matrix, then their directed sum $A\oplus B$ is the $(r+s)\times (r+s)$ matrix
\begin{equation*}A\oplus B=\left[\begin{array}{cc}A&0\\0&B\end{array}\right].\end{equation*}
The classical rational canonical form for matrices over a field can be extended to the case where the underlying field is replaced by a general division ring.
However, in the setting of division rings it is necessary to single out the center $\mathds F$ of the division ring $\mathds D$ and to disinguish the algebraic case, where the matrix satisfies a polynomial equation over $\mathds F$. Fore  details about this question the reader is referred to \cite{cohn}. We merely state the following fact which will be used later in the paper.
 If $\mathds D$ is a division ring with center $\mathds F$ and $A\in M_n(\mathds D)$, then there exists an invertible matrix $S\in M_n(\mathds D)$ such that
\begin{equation}\label{e1}SA S^{-1}=A_1\oplus\ldots\oplus A_r\oplus U,\end{equation}
where each  $A_i$, for $i=1,\ldots, r$, is an algebraic matrix over the field $\mathds F$ with a single elementary divisor $\alpha_i$, and  $U$ is transcendental matrix. In the special case where the matrix $A$ is algebraic over the field $\mathds F$ there is no transcendental part $U$ in equation (\ref{e1}). Note that, if the division ring $\mathds D$ is finite dimensional over its center $\mathds F$,  then all matrices $A\in M_n(\mathds D)$ are algebraic over $\mathds F$.  We say that a matrix $A\in M_n(\mathds D)$ is  nonderogatory or cyclic if the
 right side of  equation (\ref{e1}) has only one factor.

Following \cite[p.227]{lam} a division ring $\mathds D$ is called centrally finite if $\mathds D$ is finite-dimensional over its center $\mathds F$. Otherwise, $\mathds D$ is called centrally infinite. Also from  \cite[p.254]{lam}, if
 $\mathds D$ is a centrally finite division ring over the field $\mathds F$, with $dim_{\mathds F}\mathds D=k$, then $\mathds D\otimes_{\mathds F}\mathds D^*\cong M_k(\mathds F)$, where $\mathds D^*$ is the opposite ring of $\mathds D$. Since the tensor product $\mathds D\otimes_{\mathds F}\mathds D^*$ contains a copy of $\mathds D$ it follows that
 $\mathds D$ can be embedded in the matrix ring $M_k(\mathds F)$. This embedding can be extended in the obvious way to an embedding function
\begin{equation*}\varphi:M_n(\mathds D)\rightarrow M_{kn}(\mathds F).\end{equation*}
Therefore, the matrix ring $M_n(\mathds D)$ is isomorphic to a subalgebra $\mathds A$ of $M_{kn}(\mathds F)$. Clearly, for any matrix $A\in M_n(\mathds D)$ we have for the centralizers
\begin{equation}\label{ec}\varphi(C_{M_n(\mathds D)}(A))=C_{\mathds A}(\varphi(A))=C_{M_{kn}(\mathds F)}(\varphi(A))\cap \mathds A.\end{equation}

\

We denote the set of noncentral nilpotent and noncentral idempotent matrices in $M_n(\mathds D)$ by $\mathcal{N}_n$ and $\mathcal{E}_n$, respectively.
Note that, for any division ring  $\mathds D$, in the matrix ring $M_n(\mathds D)$ there are exactly two central idempotents, the identity matrix and the zero matrix (that is, the trivial idempotents), and only one central nilpotent, the zero matrix.

Finally, let us define the distance between a vertex  and a nonempty set of vertices of $\Gamma(M_n(\mathds D))$. Let $\mathds D$ be a division ring.
 For a nonempty set $C$ of noncentral matrices in $M_n(\mathds D)$ and a
noncentral matrix $A\in M_n(\mathds D)$ we define the distance from the matrix $A$ to the set $C$ to be $d(A, C)= min\{d(A, X) : X\in C\}$, if $A\notin C$
and $d(A, C)=0$ if $A\in C$.

\section{Main results}
It was shown in \cite[Lemma 2]{ak3} that for a field $\mathds F$ and an integer $n\geq 3$, if there is a field extension of $\mathds F$ of degree $n$ with no proper intermediate fields, then $\Gamma(M_n(\mathds F))$ is not a connected graph. In the following lemma we extend this result to centrally finite division rings.

\begin{lemma}\label{la} Let $\mathds D$ be centrally finite division ring over the field $\mathds F$, with $dim_{\mathds F}\mathds D$=k. If there exists a noncentral matrix $A\in M_n(\mathds D)$, such that the $\mathds F$-algebra generated by $A$, that is, $\langle A\rangle=\{p(A): p\in\mathds F[x]\}$  is a field extension of $\mathds F$ of degree $nk$ with no proper intermediate fields, then the commuting graph $\Gamma(M_n(\mathds D))$ is not a connected graph.\end{lemma}
{\bf Proof.} Suppose that $\langle A\rangle$ is a field extension of $\mathds F$ of degree $kn$ with no proper intermediate fields. Then, the image $\varphi(\langle A\rangle)$ of the algebra $\langle A\rangle$ by the embedding
$\varphi:M_n(\mathds D)\rightarrow M_{kn}(\mathds F)$ is a field extension of $\mathds F$ of degree $nk$ inside $M_{kn}(\mathds F).$ On the other hand, since $\varphi(\langle A\rangle)=\langle\varphi(A)\rangle$ it follows that $\varphi(A)$ is a nonderogatory matrix in
$M_{kn}(\mathds F)$. Consequently, the centralizer $C_{M_{kn}(\mathds F)}(\varphi(A))$ equals the $\mathds F$-algebra generated by $\varphi(A)$, see \cite[Corollary 1]{ak3}. Therefore, from equation (\ref{ec}) we find that
\begin{equation}\label{ec1}\varphi(C_{M_n(\mathds D)}(A))= C_{\mathds A}(\varphi(A))=\mathds F[\varphi(A)]\cap \mathds A.\end{equation}
Now, let $B$ be a noncentral matrix  in $C_{M_n(\mathds D)}(A)$, and let us find its centralizer. Since $\varphi(B)$ lies in $\in\mathds F[\varphi(A)]$ and there is no proper intermediate field between $\mathds F$ and $\mathds F[\varphi(A)]$ it follows that $\mathds F[\varphi(B)]=\mathds F[\varphi(A)]$.  That is, the $\mathds F$-algebra generated by $\varphi(B)$ is a field extension of $\mathds F$ of degree $nk$. Thus,
 the matrix $\varphi(B)$ is nonderogatory in $M_{kn}(\mathds F)$ and hence $C_{M_{nk}(\mathds F)}(\varphi(B))=\mathds F[\varphi(B)]$. So, we have
 \begin{equation}\label{ec2}\varphi(C_{M_n(\mathds D)}(B))= C_{\mathds A}(\varphi(B))=\mathds F[\varphi(B)]\cap \mathds A.\end{equation}
 Finally, since $\mathds F[\varphi(B)]=\mathds F[\varphi(A)]$, combining equation (\ref{ec1}) with equation (\ref{ec2}) it follows  that $C_{M_n(\mathds D)}(A)=C_{M_n(\mathds D)}(B)$.
This implies that the noncentral matrices in $C_{M_n(\mathds D)}(A)$ form a connected component of $\Gamma(M_n(\mathds D))$. This completes the proof. $\Box$

Our next goal is to find the distance from a noncentral matrix $A\in M_n(\mathds D)$ to the set $\mathcal{E}_n\cup\mathcal{N}_n$.

\begin{lemma}\label{lf} Let $\mathds D$ be centrally finite division ring over the field $\mathds F$ and $n\geq 3$. For a noncentral matrix $A\in M_n(\mathds D)$ if the
$\mathds F$-algebra $\langle A\rangle$ is not a field, then $d(A, \mathcal{E}_n\cup\mathcal{N}_n )\leq 1$.
\end{lemma}
{\bf Proof.} For $A\in\mathcal{E}_n\cup\mathcal{N}_n$ there is nothing to prove. So, assume that $A\notin\mathcal{E}_n\cup\mathcal{N}_n$. If the matrix $A$ is  derogatory, then by the rational canonical form we can find  a nonsingular matrix $S\in M_n(\mathds D)$ such that
\begin{equation*}A=S^{-1}(C_1\oplus\ldots\oplus C_k)S,\end{equation*}
where $k\geq 2$ and  each  $C_i$, for $i=1,\ldots, k$, is an algebraic matrix over the field $\mathds F$ with a single elementary divisor $\alpha_i$. In this case $A$ commutes with the noncentral idempotent
\begin{equation*}E=S^{-1}(I_{C_1}\oplus 0_{C_2}\oplus\ldots\oplus 0_{C_k})S,\end{equation*}
where $I_{C_1}$ is the identity matrix whose order equals the order of $C_1$ and $0_{C_i}$, for $i=2,\ldots,k$, is the zero matrix whose order equals the order of $C_i$.
Thus, $d(A, \mathcal{E}_n\cup\mathcal{N}_n )\leq 1$.

Now assume that $ A$ is non derogatory. First  suppose that the $\mathds F$-algebra  $\langle A\rangle$ generated by $A$ is a non semisimple. In this case the Jacobson radical $\mathcal{J}(\langle A\rangle)$ contains a nonzero element $J$. As is well know, see for example \cite[p.658]{lang}, the Jacobson radical of a finite dimensional algebra over a field is nilpotent. Thus, $J$ is a nilpotent matrix. Since the algebra $\langle A\rangle$ is commutative if follows that $A$ commutes with $J$ and so $d(A, \mathcal{E}_n\cup\mathcal{N}_n )\leq 1$.

Finally, assume that   $\langle A\rangle$ is semisimple and the matrix $A$ is
nonderrogatory. Since a semisimple commutative ring is isomorphic to sum of fields \cite[p.661]{lang}, if follows that
\begin{equation}\label{eq1}\langle A\rangle\cong\mathds F_1\oplus\ldots\oplus\mathds F_s,\end{equation}
where each $\mathds F_i$, for $i=1\ldots, s$, is a field.
Since $\langle A\rangle$ is not a field it follows that in (\ref{eq1}) we have $s\geq 2$. Hence,  $\mathds F_1\oplus\ldots\oplus\mathds F_s$ contains the non-trivial   idempotent $(1, 0, \ldots , 0)$. Let $\psi:\mathds F_1\oplus\ldots\oplus\mathds F_s\rightarrow \langle A\rangle $ be an isomorphism. Clearly, a non-trivial idempotent of $\mathds F_1\oplus\ldots\oplus\mathds F_s$ is mapped
by $\psi$ to a non-trivial idempotent in $\langle A\rangle$.
 As we have noted before, in the matrix ring $M_n(\mathds D)$ the central idempotent are exactly the trivial idempotent. Hence, the idempotent $\psi(1, 0, \ldots , 0)$ is not central, and the result follows. $\Box$

\begin{lemma}\label{lg} Let $\mathds D$ be a centrally finite division ring over the field $\mathds F$, with $dim_{\mathds F}\mathds D=k$ and $n\geq 3$. Assume that every field extension of $\mathds F$ of degree $nk$ has a proper intermediate field. Then,
for a noncentral matrix $A\in M_n(\mathds D)$, if the $\mathds F$-algebra $\langle A\rangle$ is a field  extension of degree    $nk$
we have $d(A, \mathcal{E}_n\cup\mathcal{N}_n )\leq 2$. If $\langle A\rangle$ is a field extension of degree degree less than $nk$ we have $d(A, \mathcal{E}_n\cup\mathcal{N}_n )\leq 1$.   \end{lemma}
 {\bf Proof.} First assume that $\langle A\rangle$ is a field extension of degree $nk$. Let $\mathds K$ be a proper intermediate field between $\mathds F$ and $\langle A\rangle$. If $B\in\mathds K-\mathds F$, then $B$ is a noncentral derogatory matrix. Hence, by the rational canonical form we can find a noncentral idempotent $E$ that commutes with $A$. Thus, $A-B-E$ is a path and so $d(A, \mathcal{E}_n\cup\mathcal{N}_n )\leq 2$.

 Now assume that $\langle A\rangle$ is a field extension of degree less than   $nk$. In this case the matrix $A$ is derogatory. Hence, it commutes with a noncentral idempotent and so $d(A, \mathcal{E}_n\cup\mathcal{N}_n )\leq 1$, and the proof in complete. $\Box$

 \

 As a consequence of  Lemma \ref{lf} and Lemma \ref{lf} we have the following result.

 \begin{Corollary}\label{c1}Let $\mathds D$ be a centrally finite division ring over the field $\mathds F$, with $dim_{\mathds F}\mathds D=k$ and $n\geq 3$. For a noncentral matrix $A\in M_n(\mathds D)$,  If the
$\mathds F$-algebra $\langle A\rangle$ is not a field or the matrix $A$ is derogatory, then $d(A, \mathcal{E}_n\cup\mathcal{N}_n )\leq 1$ \end{Corollary}

\

In  \cite{ak3} Akbari et al. established that
 for a field $\mathds F$ and an integer $n\geq 3$ the graph $\Gamma(M_n(\mathds F))$ is connected if and only if every field extension of $\mathds F$ of degree $n$ contains a proper intermediate field. On the other hand, it was established in  \cite{ak2} that if  $\Gamma(M_n(\mathds F))$ is connected, then its diameter is between four and six.
 In the following theorem, which is the main result of your paper,  we generalize this two facts to matrices over a centrally finite division ring.

\begin{theorem}\label{tc} Let $\mathds D$ be a centrally finite division ring over the field $\mathds F$, with $dim_{\mathds F}\mathds D=k$  and $n\geq 3$. The graph
$\Gamma(M_n(\mathds D))$ is connected if and only if the following condition is satisfied

$(\ast)$For any noncentral  matrix $A\in M_n(\mathds D)$ if the algebra $\langle A\rangle$ is a field extension of $\mathds F$ of degree $nk$, then there exists a proper
intermediate field between $\mathds F$ and $\langle A\rangle$.

Furthermore, if the graph $\Gamma(M_n(\mathds D))$ is connected then
\begin{equation*}4\leq diam(\Gamma(M_n(\mathds D)))\leq 6 \end{equation*}

\end{theorem}
{\bf Proof.} The only if part follows from Lemma \ref{la}. To prove the if part,  assume that condition $(\ast)$ is satisfied and let us prove that $\Gamma(M_n(\mathds D))$ is connected with diameter between four and six. First of all note that, from \cite[Lemma 2]{ak2} it follows that $diam(\Gamma(M_n(\mathds D))\geq 4$. Now, Let $A, B\in M_n(\mathds D)$ be two noncentral matrices. If condition $(\ast)$ is satisfied, then by Lemma \ref{lf} and Lemma \ref{lg} we can find noncentral matrices $X, Y\in \mathcal{N}_n\cup\mathcal{E}_n$  such that $d(A, X)\leq 2$ and
  $d(B, X)\leq 2$. Notice that, if a matrix  commutes with noncentral nilpotent matrix $N$ with index of nilpotence $\alpha$, then the matrix  also commutes with $N^{\alpha-1}$. Consequently, there is no loss of generality in assuming that that the index of nilpotence of a nilpotent matrix that commutes with a given matrix is two.

 On the other hand, following \cite[Theorem\;11]{ak2}, if $E_1, E_2\in M_n(\mathds F)$ are two noncentral idempotent matrices, then $d(E_1, E_2)\leq 2$ in the commuting graph $\Gamma(M_n(\mathds D))$. Also by  \cite[Theorem\;9]{ak2}
 if $N_1, N_2\in M_n(\mathds F)$ are two non-zero nilpotent matrices, both of index of nilpotence $2$, then $d(N_1, N_2)\leq 2$ in $\Gamma(M_n(\mathds D))$. Finally, from \cite[Proposition 1]{ga} or \cite[Lemma 3.2]{da} if $E, N\in M_n(\mathds F)$ are  such that $E$ is idempotent and $N$ is nilpotent of index of nilpotence $2$, then $d(E, N)\leq 2$ in $\Gamma(M_n(\mathds D))$. From this it follows  that $d(X, Y)\leq 2$. Hence, $d(A, B)\leq 6$, which completes the proof. $\Box$

\

Notice that, if  $n\geq3$ and any noncentral matrix $A\in M_n(\mathds D)$ satisfies $d(A, \mathcal{E}_n\cup\mathcal{N}_n )\leq 1$, then $\Gamma(M_n(\mathds D))$ is connected with diameter four. If $\mathds D$ is an algebraically closed field or a real closed field, then for any matrix $A\in M_n(\mathds D)$ the algebra $\langle A\rangle$ is not a field or the matrix $A$ is derogatory. Therefore, on applying Corollary \ref{c1} we obtain the  result already known that for an algebraically closed field or a real closed field the commuting graph is connected with diameter four.

\section{Concluding remarks} The centrally finite division algebras over a field $\mathds F$ are classified by the Brauer group $Br(\mathds F)$. For details about the Brauer group the reader
is referred to \cite{lor}. If $\mathds F$
is algebraically closed, then  the Brauer group is trivial. This is the well known fact that over  an algebraically closed field $\mathds F$, there are no finite dimensional division algebras, except $\mathds F$ itself.
If $\mathds F$ is a finite field, then $Br(\mathds F)$ is again trivial. This is a consequence of a  famous theorem of Wedderburn according to which all finite division rings are commutative. For the commuting graph of a full matrix ring over a finite field we refer the reader to \cite{da}.  If $\mathds F$ is a real closed field, then $Br(\mathds F)$ is the cyclic group of order two. In fact, there are only two centrally finite division algebras over a real closed field $\mathds F$, namely, the field $\mathds F$ itself and the quaternion algebra $\mathds H$. This is a well-know theorem of Frobenius. For an integer $n\geq 3$, let us prove that $\Gamma(M_n(\mathds H))$ is connected and find its diameter. Since for a real closed field $\mathds F$ its algebraic closure $\overline{\mathds F}$  has degree two it follows that for $n\geq 3$  the algebra $\langle A\rangle$ generated by a noncentral matrix $A\in M_n(\mathds H)$
is not a field or the matrix $A$ is derogatory.
 Hence, by Corollary \ref{c1} we conclude that $d(A, \mathcal{E}_n\cup\mathcal{N}_n )\leq 1$ and so $\Gamma(M_n(\mathds H))$ is connected with diameter four.
This result was also obtained by Akbari et al. in \cite{ak2}.

By contrast with an algebraically closed field, a real closed field  or a finite field, when $\mathds F$ is a $p$-adic field $\mathds Q_p$ the Brauer group turns out to be richer. In fact, $Br(\mathds Q_p)$ is isomorphic to  the group $\mathds Q/\mathds Z$, the  rational numbers with addition mod $1$, see \cite{serre}. Therefore, to any element $q$ of the group $\mathds Q/\mathds Z$ corresponds
a centrally finite division algebra $\mathds D_{pq}$ over the field $\mathds Q_p$.  It would be interesting to find the exact value of the diameter of the commuting graphs of matrix rings over
$\mathds D_{pq}$.

\

{\bf Acknowledgments.}

This work was supported by FCT project UID/EEA/50008/2013.

\label{}





\bibliographystyle{elsarticle-num}
\bibliography{<your-bib-database>}



\end{document}